\documentclass[a4paper,12pt]{amsart}

\usepackage{amssymb,amsbsy,amsmath,amsfonts,amssymb,amscd}
\usepackage{latexsym}
\usepackage{graphics}
\usepackage{color}
\input xy
\xyoption{all}

\newcommand\sC{{\mathcal C}}
\newcommand\sT{{\mathcal T}}
\newcommand\sD{{\mathcal D}}
\newcommand\sE{{\mathcal E}}

\newcommand\sG{{\mathcal G}}

\newcommand\sL{{\mathcal L}}

\newcommand\sX{{\mathcal X}}
\newcommand\sY{{\mathcal Y}}
\newcommand\sH{{\mathcal H}}
\newcommand\LL{{\mathbb L}}

\newcommand\Ga{\Gamma}
\newcommand\De{\Delta}
\newcommand\ga{\gamma}
\newcommand\de{\delta}

\DeclareMathOperator{\Ext}{Ext}

\newcommand{\CC}{\ensuremath{\mathbb{C}}}

\newcommand{\ZZ}{\ensuremath{\mathbb{Z}}}

\newcommand{\hol}{\ensuremath{\mathcal{O}}}

\newcommand{\PP}{\ensuremath{\mathbb{P}}}

\newcommand{\ra}{\ensuremath{\rightarrow}}

\def\eea{\end{eqnarray*}}
\def\bea{\begin{eqnarray*}}

\newcommand\dual{\mathrel{\raise3pt\hbox{$\underline{\mathrm{\thinspace d
\thinspace}}$}}}
\newcommand\qe{\ifhmode\unskip\nobreak\fi\quad $\Box$}       

\def\BOX{\hfill\lower.5\baselineskip\hbox{$\Box$}}

\newtheorem{theo}{Theorem}[section]
\newtheorem{remarkk}[theo]{Remark}
\newenvironment{rem}{\begin{remarkk}\rm}{\end{remarkk}}

\newtheorem{defin}[theo]{Definition}

\newtheorem{prop}[theo] {Proposition}

\newtheorem{lemma}[theo]{Lemma}
\newtheorem{example}[theo]{Example}

\newcommand{\Proof}{{\it Proof. }}
\begin{document}

\title[Cyclic curve covers]{ Irreducibility of the space of cyclic
covers of algebraic curves of fixed numerical
type  and the  irreducible components of
$Sing (\overline{\mathfrak M_g})$}
\author{ F. Catanese}

\thanks{The present work took place in the realm of the DFG
Forschergruppe 790 "Classification of algebraic
surfaces and compact complex manifolds". }

\date{\today}

\maketitle

{\em Dedicated to Shing Tung Yau with friendship and admiration on
the occasion of his 60-th birthday.}

\section*{Introduction}

The first main purpose of this  article is to prove  irreducibility
(equivalently,  connectedness) results
(Theorems \ref{irreducible}, \ref{stableirred}) for the space of cyclic covers
   of a fixed numerical type between  complex
projective curves, first in the case of smooth curves, and then in the
case of stable curves.

In the smooth case the numerical type of a cyclic cover $ C \ra C'$  is
given  by the order $d$ of the cyclic group
$G$, by the genus $g$ of the covering curve $C$, and by the branching
datum, i.e., the equivalence class, for
the natural action of $Aut (G) = (\ZZ/d)^*$, of the branching sequence
$(k_1, \dots , k_{d-1})$, where
$k_i$ is the  number  of branch points on $C'$ such that their local
monodromy is the  element $i$ of the
group $G \cong \ZZ / d$.

This result, in the case where $C$ is smooth and $d$ is prime, was
obtained long ago by Cornalba (\cite{cor1});
Barbara Fantechi \footnote{Talk at MSRI, February 2009.} observed
that Cornalba's proof does partially extend for $C$  smooth, 
but, for arbitrary  $d$, not for any numerical type.

Similar  results were obtained   by Biggers and Fried in
\cite{bf}, but  concerning only the case where $ C \ra C'$
is unramified and cyclic.

Observe that, much more generally, given any finite group $G$,
the space of Galois covers $C \ra C'$ with group $G$, with $C$ smooth, 
 and with a fixed topological type
is connected. This follows from Teichm\"uller theory, as shown in
Proposition 4.13  of  \cite{FabIso} (see also sections 5 and 6 of \cite{cime}
for related topics). 

This result reduces the question of determining the possible
topological types to a question in group theory, namely
finding the equivalence classes of monodromies $\mu : \pi_h :=  \pi_1 (C') \ra G$,
for the action of $Aut (G)$ and for the action on the source of the Teichm\"uller group $\sT_h = Out
(\pi_h)$. In general the above numerical invariants of the cover 
(the branching sequence is labeled here by the conjugacy classes in $G$) do not determine the topological type,
but our first theorem says that they do determine it in the case where $G$ is cyclic.

The second main purpose of this article is on the other hand  (Theorem \ref{stablesing}) the
description of  an irredundant irreducible
decomposition of
$Sing (\overline{  \mathfrak M_g })$.

The determination  of an irredundant irreducible
decomposition for $Sing(\mathfrak M_g )$ was  obtained by
Cornalba in \cite{cor1} and \cite{cor2};
we give here a slightly shorter proof of Cornalba's description (\ref{cornalba}).

Our main novel contribution  is  then the   complete explicit determination
of the irreducible components
which are contained in the
boundary
$\overline{  \mathfrak M_g } \setminus  \mathfrak M_g $.

Here a brief description of  the contents of the paper.

In the first section we recall the general theory of cyclic coverings
between normal varieties (and with
factorial base), expecially important to construct  families of
cyclic covers: here a more general theory of
abelian covers and their invariants was developed by Pardini in
\cite{pardini}, extending the case of
bidouble covers developed in \cite{ms}.

In the second section we prove the irreducibility result for the
space of cyclic covers of smooth curves with a
fixed topological type (Theorem \ref{irreducible}) and for arbitrary
degree $d$.

In the third section we describe  the relation of the above result
with the determination of the
irreducible components of the  singular locus of the moduli space of curves
$\mathfrak M_g$.

In fact, for $ g \geq 3$, with exception of the locus of 
hyperelliptic curves of genus
  $3$, all the loci that  one
obtains by letting $d$ be a prime number are closed subvarieties of $
\mathfrak M_g$ of codimension $ c
\geq 2$, hence they are irreducible sets whose union is $Sing ( 
\mathfrak M_g )$.
We give a shorter proof of Cornalba's theorem, using (in Proposition
\ref{normalizer}) numerical inequalities
instead  of degeneration arguments.

We restate, with some additions (which are then needed in section 4), Cornalba's main result in Theorem
\ref{cornalba}.

In the fourth  section we consider  the pair of a stable curve $C$
and an  automorphism $\ga$, and we
describe some numerical (and combinatorial) invariants of the pair $(C, \ga)$. A
fundamental difference here, when trying  to
describe irreducible components of the `space' of such coverings $ C
\ra C /\ga$, is that the topological
type is not constant on connected components, but each such component
is a finite union of locally closed
irreducible sets (corresponding to a given topological type),
only one of them being open and dense.

We call the pairs $(C, \ga)$ in this dense open set {\bf maximal},
and we give a
complete numerical-combinatorial description of such maximal coverings in the
simpler case where the order $d$ of the
automorphism
$\ga$ is a prime number (Theorem \ref{stableirred} deals only with
the case where $d$ is
a prime number, we do not treat here the general case).

We also describe a
`simplification' algorithm showing how to obtain, from a fixed
numerical-combinatorial type, the associated maximal
numerical-combinatorial type.

Theorem \ref{stableirred} allows us to  completely  describe  in the
final section (Theorem \ref{stablesing}) the
irredundant irreducible decomposition of $Sing (\overline{  \mathfrak M_g })$.

The irreducible components of $Sing (\overline{  \mathfrak M_g })$
fall into two types: the ones of the first
type   are the closures of the irreducible components of
$Sing (\mathfrak M_g )$, while the ones of the second type are the
    components  which are completely contained in the
boundary
$\overline{  \mathfrak M_g } \setminus  \mathfrak M_g $.

The determination of these latter components is  different than in the case
of  $Sing (\mathfrak M_g )$,
and we have the following phenomenon.

The stable curves $C$ with an elliptic tail $E$ (i.e., $E$ is a smooth
elliptic curve intersecting
$C'' : = \overline {C \setminus E}$ in one point $P$)  form a divisor,
and admit an involution
(an automorphism
of order $2$) which is the identity on $C''$, and on $E$ is
multiplication by $-1$ when one chooses $P$ as
the origin of the elliptic curve. If a curve $C$ has several elliptic
tails, the corresponding involutions
are in the centre of $Aut (C)$ and do not contribute to $Sing
(\overline{  \mathfrak M_g })$.

It turns out that the irreducible components  which are completely
contained in the
boundary correspond  to cyclic automorphisms of prime order which are the identity on
all components except one,
 are of maximal numerical-combinatorial type,  and, if the order $d$ equals 2, they
act as the identity on the elliptic tails. For two exceptional cases, which
we explicitly describe, we do not get irreducible components, but
proper  subsets.

In a brief final section we mention some  problems and results related
to the automorphism group $ Aut (C)$ of a stable curve $C$ of genus $g$
(cf. \cite{VO1},\cite{VO2},\cite{VO3}) .

\section{Cyclic covers of factorial varieties}

Let $X, Y$ be  normal complex projective varieties, and let $\CC(X)$,
$\CC (Y)$ be their respective  function
fields.

Assume that  $\CC(X)$ is a cyclic  Galois extension of $\CC (Y)$, and
denote by $$G \cong \mu_d : = \{
\zeta \in \CC | \zeta^d = 1\}$$ its Galois group,  by $\ZZ/d$ the
group of characters $$ \ZZ/d \cong \{ \chi
|\exists m \in \ZZ/d ,  \chi (\zeta) = \zeta^m\}.$$

For each character $\chi$ of order  $d$ the extension is given by
$$\CC(X) =  \CC (Y)(w) , w^d = f(y) \in \CC (Y), $$ where $w$ is a
$\chi$-eigenvector.

Assume now that $Y$ is  factorial, so that $f$ admits a unique prime
factorization as a fraction of pairwise
relatively prime sections of line bundles, and we can write
$$  w^d = \frac{\prod_i \sigma_i^{n_i} } {\prod_i \tau_j^{m_j} } .$$

Write now $$ n_i =  N_i + d n'_i , \ \  m_j = - M_j + d m'_j $$ with
$0 \leq N_i, M_j \leq d-1$ and set
$$ z : = w \cdot  \prod_i \sigma_i^{- n'_i} \prod_i \tau_j^{m'_j}.$$

Whence
$z$ is a rational section of a line bundle on $Y$ and we have

$$   z^d = \prod_i \sigma_i^{N_i} \prod_i \tau_j^{M_j}  . $$

We put together  the prime factors which appear with the same
exponent and write:
$$   z^d = \prod_{i=1}^{d-1} \delta_i^{i}  . $$

Here each factor $\delta_j$ is reduced, but not irreducible, and
corresponds to a Cartier divisor that we
shall denote $D_j$. The local monodromy around $D_j$ is easily
calculated since, if we take a small loop
$$ \de_j  = e^{ 2 \pi \sqrt{-1} \ \theta} , \theta \in [ 0, 1]$$ it 
lifts to the path
    $$z= z_0   \cdot  e^{ \frac{2 \pi j \sqrt{-1} \  \theta}{d}}$$

and its monodromy is
$$ z_0  \mapsto e^{ \frac{2 \pi  j \sqrt{-1} }{d}} z_0.$$

This shows that $D_i$ is exactly the divisorial part of the branch
locus $ D : =  \sum_i  D_i $ where the local
monodromy is the $i$-th power of the standard generator
$\ga : = e^{ \frac{2 \pi   \sqrt{-1} }{d}} $ of $G \cong \mu_d$.

In the following we shall write characters additively, in the sense
that we view them as
$ G^{\star}  = Hom (G, \ZZ / d)$. We notice also that what said above
applies to any character $\chi$: to
$\chi$ we associate  the normal covering
$$ Z_{\chi} : = X / ker (\chi). $$

We have then a linear equivalence
   $$ (*)  \ d L _{\chi} \equiv \sum_i  \overline{\chi(i)} D_i $$
   where $\overline{m}$, for $m \in \ZZ / d$, is the unique 
representative of the residue class lying
   in $ \{ 0,1, \dots , d-1\}$.

   We observe for further use the following formula:
   $$  (I) \ \   \overline{\chi(i)}  +  \overline{\chi'(i)} =
\overline{(\chi + \chi')(i)} + \epsilon^i_{\chi, \chi'}
d,$$
   where $ \epsilon^i_{\chi, \chi'}  \in \{0,1\}$.

The following theorem is in part a special case of the structure
theorem for Abelian coverings
   due to  Pardini ( \cite{pardini}, the existence  result 
in terms of the basic linear equivalences was however
already obtained by Comessatti in
\cite{comessatti}): but here  is  explicitly stated the
irreducibility criterion for the covering.

\begin{theo}\label{cyclic} i) Given a factorial variety $Y$, the
datum of a pair $(X, \gamma )$ where $X$ is
a normal variety and $\ga $ is an automorphism of order $d$ such
that, $G$ being the subgroup generated
by $\ga$, one has  $X / G \cong Y$, is equivalent to the datum of
reduced effective divisors $D_1, \dots
D_{d-1}$ without common components, and of a divisor class $L$ such
that we have the following linear
equivalence
$$ (*)  \ d L \equiv \sum_i i D_i $$ and moreover, setting $ m : =
G.C.D. \{ i | D_i \neq 0\} $, either

(**) $ m=1$ or, setting $ d = mn$, the divisor class
$$ (***)  \ L' : = \frac{d}{m} L -  \sum_i \frac{i}{m}  D_i $$ has
order precisely $m$ in the Picard group.

ii) In fact, if $\LL$ is the geometric line bundle whose sheaf of
regular sections is $\hol_Y(L)$, then
$X$ is the normalization of the singular covering
$$ X' \subset  \LL, X' : = \{ (y,z)|  z^d = \prod_{i=1}^{d-1}
\delta_i^{i} (y) \}. $$ And $\ga$ acts by $ z
\mapsto e^{ \frac{2 \pi   \sqrt{-1} }{d}} z $.

iii) The scheme structure of $X$ is explicitly given as
$$ X : = Spec ( \hol_Y \oplus (\bigoplus_{\chi \in G^{\star}
\setminus \{0\}} \hol_Y (- L_{\chi}) ))$$ where
the divisor classes $L_{\chi}$ are recursively determined by $L_1 : =
L$, and by $ L_{\chi + \xi} \equiv
L_{\chi} + L_{\xi} + \sum_i  \epsilon^i_{\chi, \xi} D_i$.

And where  the ring structure is given by  the multiplication maps
$$\hol_Y (- L_{\chi}) \times \hol_Y (-
L_{\xi}) \ra \hol_Y (- L_{\chi+ \xi})$$ determined by the section
$$\prod_i \delta_i^{ \epsilon^i_{\chi,
\xi}} \in  H^0 ( \hol_Y (- L_{\chi+ \xi} + L_{\chi} + L_{\xi})).$$

\end{theo}

\begin{rem} We can write more suggestively the ring structure as
$$z_{\chi} \cdot z_{\xi}   =  z_{\chi+ \xi} \prod_i \delta_i^{
\epsilon^i_{\chi, \xi}},$$ where $z_{\chi} $is
thought as a fibre variable on the geometric line bundle
$\LL_{\chi}$.

In other words, $X$ is thus embedded in $\bigoplus_{\chi \in
G^{\star} \setminus \{0\}} \LL_{\chi}$.

\end{rem}

\Proof We need only to show that $X$ is irreducible iff either (**)
or (***) holds.

But  a covering is irreducible if and only if the covering monodromy
is transitive, and this clearly holds if
$m=1$.

Else, we consider the quotient of $X$ by the subgroup $H$ of $G$
generated by the inertia subgroups (i.e., $
H = m \ZZ / d \ZZ$). Set $ Z_m = X /H$, so that (**) holds,  hence
$X$ is irreducible iff $Z_m$ is irreducible.

More concretely, if  $ m > 1$,  set $ u : = \frac{z^n}{\prod_i
\delta_i^{\frac{i}{m}}}$, so that we  have a
factorization of the covering given by
$$u^m = 1, \  z^n  = u \prod_i \delta_i^{\frac{i}{m}}.$$

This means that we take $ X \ra Z_{m} \ra Y$, and the last covering is
\'etale with group $\cong \ZZ / m$.

The last covering is irreducible then if and only if the divisor
class $L'$ has order exactly $m$.

\qed

\begin{rem}\label {families} The explicit description of cyclic
covers (resp. : abelian covers) allows to
construct families of varieties with a cyclic automorphism.

Let $\sY \ra \sT$ be a proper morphism with projective fibres, such
that $\sY$ is factorial, and consider
relative effective Cartier  divisors $\sD_i$ and a relative Cartier
divisor $\sL$ such that
$$  (*)  \ d \sL \equiv \sum_i i \sD_i .$$

Then we have a finite Galois morphism $\sX \ra \sY$ over $\sY \ra
\sT$ with an automorphism $\ga$ over
$\sY \ra \sT$ such that $\ga$ generates the Galois group $G$ and $\sX
/ G =  \sY$. In particular, for each
fibre $X_t$,
$X_t / G = Y_t$.
\end{rem}

\section{Cyclic covers of curves and their invariants}

In this section $C$ will be a projective complex curve of genus $g$,
and $\ga$ an automorphism of $C$ of
order exactly $d$.

In this situation we can consider some obvious numerical invariants.

\begin{defin} Let  $G \cong \ZZ / d $ be the subgroup generated by
$\ga$, and set $C' : =  C / G$,
$h: = $ genus$(C')$.

Denote by $k_i : = deg (D_i)$ for $ i=1, \dots, d-1$, and by $(k_1,
\dots k_{d-1})$ the branching sequence
of $\ga$.

\end{defin}

Observe  that if we change the chosen generator for $G$, then the
covering does not change, only the
identification of the Galois group with $\ZZ / d$ changes. Then we
get another  branching sequence,
obtained by multiplying the indices of the elements of the given
sequence with a fixed element  $r \in (\ZZ
/ d)^*$ (for example, for $r=-1$, we get the new sequence $( k_{d-1},
\dots, k_1))$.

\begin{defin} We shall call a branching datum an equivalence class of
branching sequences for the above
action of $(\ZZ / d)^*$ by  index multiplication. We shall denote it
by $[(k_1, \dots k_{d-1})]$.
\end{defin}

The next question is: which branching data do actually occur for $g,d$ given?

   A first  restriction is  given by the Hurwitz formula
$$ 2(g-1) = d  \{  2 (h-1) +  \sum_i  k_i  (1 - \frac{ GCD(i,d)}{d})
\},$$ which determines the genus of the
quotient curve in terms of $g,d$ and of the branching datum: it is
necessary that  $h$ be a non negative
integer number.

Moreover, by the previous theorem \ref{cyclic} a necessary and
sufficient condition is, once the above
formula yields $h \geq 0$, that
$$ (*)   \sum_i  k_i   i \equiv 0 \in \ZZ / d,$$ since in the Picard
group of a curve a divisor is divisible by
$d$ iff its degree is divisible by $d$.

This motivates the following

\begin{defin} A  branching datum corresponding to a sequence  $[(k_1, 
\dots k_{d-1})]$
is said to be admissible for $d$ and $g$ if (*) holds, and moreover
$$     h : =  1 +  \frac{2(g-1)}{2d}  -  \frac{1}{2} \sum_i  k_i  (1 
- \frac{ GCD(i,d)}{d}) $$
is a non negative integer.
\end{defin}

Unless otherwise specified, we shall  consider, given integers $d,g$, 
only admissible
branching sequences.

One can view the same result from the point of view of Riemann's
existence theorem: such pairs $(C, \ga)$
are determined by  the following data: a curve $C'$ of genus $h$,
divisors $D_i$ for all $i  \in \ZZ / d$ and
a surjective homomorphism $\psi : H_1 (C' \setminus D, \ZZ) \ra  \ZZ
/ d$ such that the image of a small
circle around a point $ p \in D_i $ maps to the class of $i$ in $ \ZZ / d.$

We have the following exact sequence, where we write $D_1 = p_1 +
\dots + p_{k_1}$, $D_2 =  p_{k_1 + 1}
+ \dots p_{k_1  + k_2} , \dots $,
$$ (**) 0 \ra  A:=  (\oplus _j \ZZ p_j ) / \ZZ (\sum_j p_j ) \ra  H_1
(C' \setminus D, \ZZ) \ra H_1 (C' , \ZZ)
\cong \ZZ^{2h} \ra 0$$ which admits several splittings.

The condition $(*)$ pertains to the relation holding in the subgroup
$A$, while we may choose a splitting
such that
$\ZZ^{2h}$ maps onto $\ZZ/ d$.  Topologically, this means that we
choose a special symplectic basis of $
H_1 (C' , \ZZ)$, such that the points $p_j$ lie in the complement  of
the corresponding canonical dissection
of the curve $C'$, and then we take a disk $\De$ contained in this
complement  and containing the branch
divisor $D$. Hence the ramified covering is just obtained glueing
together a ramified covering of $\De$
with an unramified covering of $C' \setminus \De$.

For later use, we observe that, if $d$ is a prime number $p$, then
the Hurwitz formula is easier to write and
we have

$$ 2(g-1) = p  \{  2 (h-1) +  k  (1 - \frac{ 1}{p}) \}, \ k : = \sum_i  k_i.$$

\begin{theo}\label{irreducible} The pairs $(C,G)$ where $C$ is a
complex projective curve of genus $g \geq
2$, and $G$ is a finite cyclic group of order
$d$ acting faithfully on $C$ with a given branching datum $[(k_1,
\dots k_{d-1})]$ are  parametrized by a
connected complex manifold  $\sT_{g;d,[(k_1, \dots k_{d-1})]}$ of
dimension $3 (h-1) + k $,  where $k : =
\sum_i  k_i .$

The image  $\mathfrak M_{g;d,[(k_1, \dots k_{d-1})]}$ of
$\sT_{g;d,[(k_1, \dots k_{d-1})]}$ inside the
moduli space
   $\mathfrak M_g$ is a closed subset of the same dimension $3 (h-1) + k $.

\end{theo}

\Proof Consider the Teichm\"uller space of curves $C'$ of genus $h$
with $k$ marked points, and take a
homomorphism
$\psi$ of the first homology group of  $C'$ onto $\ZZ/ d$ sending the
generator $p_j$ of the  subgroup
$\ZZ p_j$, for $ k_1 + \dots + k_{i-1} + 1 \leq j  \leq  k_1 + \dots
+ k_{i} $ to $i \in \ZZ / d$. Since the
topological type of $C'$ is fixed for the Teichm\"uller family, we
choose a fixed splitting of $(**)$ and
consider the surjection   $ H_1 (C' , \ZZ) \ra \ZZ / d$ corresponding
to a fixed primitive element $\Psi \in
H^1 (C' , \ZZ)$.

Recall that the symplectic group $Sp (2h, \ZZ)$ acts transitively on
the set of such primitive elements.

Applying remark \ref{families} to our situation, we have a family of
cyclic covers of the curves $C'$,
parametrized by Teichm\"uller space $\sT$.

We want now  to show that every  pair $(C,G)$ as in the statement
occurs in our family. To this purpose,
denote by $C'$ the quotient curve, choose an isomorphism of $G$ with
$\ZZ / d$ and a marking of the
branch points so that the branching sequence is $(k_1, \dots ,
k_{d-1})$, and the divisor $D_i$ consists of
the  points $p_j$ with $ k_1 + \dots + k_{i-1} + 1 \leq j \leq  k_1 +
\dots + k_{i} $.

We want to show that for a suitable diffeomorphism of $C'$ which
leaves the disk $\De$ pointwise  fixed  we
can transform the resulting homomorphism  $\Psi$ into the standard
one we have chosen.

To this purpose we take a product of Dehn twists over loops supported
in $C' \setminus \De$, and we
observe that these generate  the mapping class group of $C'$. Since
the mapping class group maps onto the
symplectic group, our first statement is thus proven.

The group $G$ has a linear representation on the tricanonical vector
space $H^0 ( \hol_C ( 3 K_C))$, and on
the family  $\sT_{g;d,[(k_1, \dots k_{d-1})]}$ the dimension of the
eigenspaces is semicontinuous, hence
constant, since the sum is the  fixed integer $5g-5$.

We consider a vector space $V$ of dimension $5g-5$ with a  linear
action of $G$ having eigenspaces of the
given dimensions. Inside $\PP (V^{\vee})$ we look at the locally
closed subset $\sH_g$ of the Hilbert
scheme corresponding to tricanonically embedded smooth curves of
genus $g$. Inside  $\sH_g$ we consider
the closed subset  $\sH_g^G$ of subschemes which are $G$-invariant.

We claim that the image $W$ of $\PP GL(5g-5) \times \sT_{g;d,[(k_1,
\dots k_{d-1})]}$ is closed inside
$\sH_g^G$. If $t_0$   is in the closure of $W$, we can choose an
analytic curve map  $(T, t_0) \ra
\sH_g^G$ such that $T$ is biholomorphic to a 1-dimensional disk and
there is a point $t_1$  such that
$t_1$ maps to $W$.

Over $T$ we have a family $\sC$ of curves of genus $g$ with an
automorphism $\ga$ of order $d$, where
$\ga$ generates $G$. Then $\sC / G : = \sC'$ is a family of curves of
genus $h$, and  each point $p$ in a
fibre
$C'_{t_2}$ belongs to a section $\sigma_p$ of $\sC' \ra T$, for which
the stabilizer (i.e., the associated local monodromy subgroup)  of $\sigma_p \cap
C'_{t}$  is, $\forall t$, equal to the stabilizer of $p$ (this follows
since the action of $G$ on  $\sC$ may locally be
linearized, and the action on the  base is the identity) and moreover the
character of the tangent  representation is
also the same.

Hence the number of points $k_i$ is constant for each $ i = 1, \dots
, d-1$, and we have proven that $W$ is
closed, hence the image  $\mathfrak M_{g;d,[(k_1, \dots k_{d-1})]}$
of $\sT_{g;d,[(k_1, \dots k_{d-1})]}$
inside the moduli space
   $\mathfrak M_g$ is also closed.

   Now the space $Kur(C)^G$ of  $G$ invariant deformations of our
curves $C$ correspond to the
submanifold of the Kuranishi
   family of $C$ obtained considering the subspace $H^1(\Theta_C)^G$ of
$H^1(\Theta_C)$.

   As shown in Pardini's article, there is an isomorphism between
$H^1(\Theta_C)^G$ and
   $H^1(\Theta_{C'}(-D))$. This shows that the map between
$\sT_{g;d,[(k_1, \dots k_{d-1})]}$ and
$Kur(C)^G$
   is a local biholomorphism, hence our  assertion on the dimension of
$\mathfrak M_{g;d,[(k_1, \dots
k_{d-1})]}$.

\qed

\begin{rem}\label{numbers}

The general curve $C$
inside  $\mathfrak M_{g;d,[(k_1, \dots
k_{d-1})]}$ has $G$ as a maximal cyclic group of automorphisms unless possibly when
we are in the cases
\begin{enumerate}
\item
$h=2$, $k=0$;
\item
$h=1$, $k=2$;

\item
$h=0$, $k=3$ or $4$.

\end{enumerate}

We shall later see in theorem \ref{cornalba} that  the only 
occurring case is the third.
\end{rem}

\Proof
 If the group $G$ is not a maximal cyclic subgroup of
automorphisms, there would exist a
   nontrivial automorphism of $C'$ leaving the divisor $D$ invariant.

  Since the divisor $D$ can be chosen freely, this is a contradiction if $C'$ has genus $ h \geq 3$,
  and also in the case where $h=2$
and  $D$ is non trivial.

  By our assumption,
in the case where $C'$ has genus $h=1$, $D$ contains at least two points,
and, if there are at least three points, we use that 
 any automorphism of finite order on a general elliptic curve has order 2,
hence $\phi$ leaves invariant a proper   subset $D' \subset D$,
with the property  that  the group $H$ of automorphisms leaving $D'$ invariant is finite.
Then we derive a contradiction choosing the other points of $D$ not 
to build a union of $h$-orbits for any element
$h \in H$. 

The same contradiction is derived in the case where $h=0$ and $k \geq 5$.

If there is a nontrivial automorphism $\phi$ sending $D$ to itself,
and $\phi$ has only one orbit on $D$, then $\phi$ is cyclic of order $k$
and the cross ratios of the $k$ points satisfy an algebraic relation.

Otherwise there are at least two orbits, and a proper  invariant subset $D' \subset D$,
with at least 3 elements. Since the group $H$ of automorphisms leaving $D'$ invariant is finite,
again we can choose the other points of $D$ not to build a union of $h$-orbits for any element
$h \in H$. 

\qed

\section{Irreducible components of $Sing(\mathfrak M_g)$.}

We begin this section with the obvious but very important observation
that if a curve $C$ of genus $g \geq
2 $ has a nontrivial automorphism group, then (since $Aut(C)$ is
finite) it has a non trivial automorphism
$\ga$ of finite order $d$, and in particular it has an automorphism
of prime order $ p >1$  (write $ d = p
\cdot m$).

By Theorem \ref{irreducible} and the ensueing remark the locus of curves with an automorphism
of the same type as $\ga^m$ is
strictly bigger than the locus of curves with an automorphism of the
same type as $\ga$, unless we are in 3 a priori possible special cases.

Therefore the locus of curves in $\mathfrak M_g$ with nontrivial
automorphisms is the union of the
irreducible closed subsets
$\mathfrak M_{g;p,[(k_1, \dots k_{p-1})]}$ (where $p$ is a prime number
and the sequence $(k_1, \dots k_{p-1})$ is $(g,p)$ admissible).

\begin{rem}{\bf(Codimension of these loci) }

$\mathfrak M_{g;p,[(k_1, \dots k_{p-1})]}$ has dimension $ 3 (h-1) +
k$, while Hurwitz' formula reads out
as
$$  2 (g-1) =   2p (h-1) + k (p-1) \Leftrightarrow  3 (g-1) =   3p
(h-1) +  \frac{3}{2} k (p-1),$$ thus the
codimension of this locus inside $\mathfrak M_g$ equals $$ c : =
3(g-1) - 3 (h-1) -  k =  3 (p-1) (h-1) + k \{
\frac{3}{2} (p-1) - 1) =$$
$$ =   3 (p-1) (h-1) + k(p-1) + k  \{  \frac{1}{2} (p-1) - 1).$$

If $h \geq 2$, then $ c \geq  3$, if $h=1$ we get $ c \geq 4$ unless
$p=2$, in which case $ k = 2 (g-1)$, whence $
c = (g-1)$, and $ c \geq 2$ for $g \geq 3$.  If $h =0$, $ c =  (k-3)
(p-1) + k  \{  \frac{1}{2} (p-1) - 1)$, which
is,  for $ p \geq 3$, $\geq 2$ unless
$ k=3, p = 3$: but then $g=1$.

Finally, if $h=0$, $p=2$, $ c = \frac{k}{2}-3 \geq 2 $ unless   $k=6$ or $k=8$.

Hence the only exceptions to $ c\geq 2$ are:  all
curves of genus  $2$ are hyperelliptic,
hyperelliptic curves of genus $3$ form  a divisor inside $\mathfrak M_3$,
and double covers of elliptic curves form  a divisor inside $\mathfrak M_2$.
\end{rem}

The next question is whether writing  the locus of curves in
$\mathfrak M_g$ with nontrivial
automorphisms as the union of the irreducible closed subsets
$\mathfrak M_{g;p,[(k_1, \dots k_{p-1})]}$
   is an irredundant  irreducible decomposition.

\begin{prop}\label{normalizer}
Assume that $p,q$ are prime numbers, and that the component
$\mathfrak M_{g;p,[(k_1, \dots k_{p-1})]}$
is contained in another (different) component $\mathfrak M_{g;q,[(k'_1, \dots k'_{q-1})]}$.

For general $C \in \mathfrak M_{g;p,[(k_1, \dots k_{p-1})]}$,  set
$A:= Aut (C)$, denote by $\ga \in A$
the automorphism of order $p$ with the first topological type, and
by $\ga' \in A$
the automorphism of order $q$ with the second topological type.
Set $ C_1 : = C / \ga,  C_2 : = C / \ga',  C'' : = C / A$, and let
$h_j$ be the genus of $C_j$,
while we let $h'' $ be the genus of $C''$; set finally $k : = \sum_i k_i,
k' : = \sum_j k'_j$.
Let $G$ be the cyclic group generated by $\ga$.

Then the normalizer $A'$ of $G$ in $A$ is strictly bigger than $G$.
\end{prop}

\Proof
Denote by $ a$ the cardinality of $A$. Our assertion holds trivially
if $ a = 2p$,
hence we shall assume that $ a \geq 3p$.

If an element $\alpha \in A$ stabilizes a point which is fixed by
$\ga$, then $\alpha , \ga$
generate a cyclic subgroup and, by our previous remark, $\alpha \in G$,
with the possible exceptions  $h_1 = 1, k=2; h_1 = 0, k=3 \ {\rm or} \ 4$.

If in these exceptional cases $\alpha \notin G$, then obviously our assertion holds,
hence we may reduce to consider the case where such an element $\alpha$ necessarily belongs to $G$.

In this case $C \ra C_1$ is branched in $k$ points, while the covering
$C_1\ra C''$ is
unramified over the images of these points.


Assume then that the normalizer of $G$ in $A$ is equal to $G$.

Then the $k$ branch points of $C \ra C_1$ have different images in $C''$,
hence, if $k''$ denotes the number of branch points of $ C \ra C''$,
$k'' \geq k$.

Hurwitz' formula yields then:
$$ 2 (g-1) = p [2 (h_1-1) + k \frac{(p-1)}{p}] \geq  a  [2 (h'' -1)
+   k \frac{(p-1)}{p} +
\frac{1}{2}(k'' -k)],$$ whence
$$ (***) \ \  2 (h_1-1) + k \frac{(p-1)}{p} \geq  \frac{a}{p} [ 2 (h'' -1)  +
k \frac{(p-1)}{p} + \frac{1}{2}(k''-k)].$$

By counting the number of moduli for the two families, we get that
$$ 3 (h_1-1) + k \leq  3 (h'' -1) + k''  \Leftrightarrow  2 (h_1-1)
\leq  2 (h'' -1) + \frac{2}{3} (k''-k).$$
The previous inequality $(***) $ is contradicted if $h'' \geq 1$, since
$\frac{a}{p} \geq 3.$

If $h'' = 0$, then by the above inequality $ k''-k \geq 3 h_1$,
hence
$$  2 (h_1-1)  \geq  \frac{a}{p} [ -2 + \frac{3}{2} h_1]   +
(\frac{a}{p} - 1) k \frac{(p-1)}{p} =$$
$$   \frac{3a}{2p}  [h_1 - 4/3]   +  (\frac{a}{p} - 1) k \frac{(p-1)}{p} .$$

If $h_1 \geq 2$ we get, since   $\frac{a}{p} \geq 3$,  $ h_1 \leq \frac{8}{5}$, a contradiction.
 
If $h_1 = 1$ we get
$$    \frac{1}{2}\frac{a}{p}   \geq  (\frac{a}{p} - 1) k \frac{(p-1)}{p} ,$$
and, since $k \geq 2$, we obtain $ \frac{1}{2} \frac{a}{p}   \geq 
(\frac{a}{p} - 1)$,
contradicting $\frac{a}{p} \geq 3$.

If instead $h_1 = 0$ we get
$$  2 (\frac{a}{p} - 1)  \geq (\frac{a}{p} - 1) k \frac{(p-1)}{p}
\Leftrightarrow 0  \geq
-2 + k \frac{(p-1)}{p},$$
hence $ g \leq 1$, a contradiction.

\qed

In the following theorem, which is the main result of Cornalba in
   \cite{cor1} and \cite{cor2}, we shall use the same notation
introduced in Proposition \ref{normalizer}.

\begin{rem}
In order to understand some statement in the following theorem, 
observe that, for $g \geq 3$,
the singular locus of  $\mathfrak M_g$ consists of the union of all 
the closed irreducible sets
$\mathfrak M_{g;p,[(k_1, \dots k_{p-1})]}$, with exception of the divisor
$\mathfrak M_{3;2,[(6)]} \subset \mathfrak M_3$, corresponding to the 
locus of hyperelliptic curves.

This holds since $\mathfrak M_g$ is locally the quotient of the 
Kuranishi family of a curve $C$
by the action of $ Aut(C)$. And, by Chevalley's theorem 
(\cite{chevalley}), this quotient is smooth if and only if
the action of $ Aut(C)$ is generated by pseudoreflections.  A 
pseudoreflection has a fixed locus which is a divisor,
hence the only pseudoreflection which occurs is the hyperelliptic 
involution in genus $g=3$.

An irreducible closed set $\mathfrak M_{g;p,[(k_1, \dots k_{p-1})]}$ 
is then an irreducible
component of  $ Sing (\mathfrak M_g )$ if it is not properly 
contained in another irreducible
set $\mathfrak M_{g;q,[(k'_1, \dots k'_{q-1})]}$, which is different 
from $\mathfrak M_{3;2,[(6)]}$.

\end{rem}

\begin{theo}\label{cornalba}
Assume that $p,q$ are prime numbers, and that the component
$\mathfrak M_{g;p,[(k_1, \dots k_{p-1})]}$
is contained in a different  component $\mathfrak M_{g;q,[(k'_1, \dots k'_{q-1})]}$.

Let $A' \neq G$  be the normalizer of $G$ in $A$, and set $G' : = A' /G$.

We have exactly the following cases:

\begin{enumerate}
\item
$h_1 = 2, k = 0$ (hence $g=p+1$) and $G' \cong \ZZ/2$ is generated by
the hyperelliptic involution; $A'$ is a dihedral group $D_p$.

For $p=2$ we have $g=3$
and $A' \cong (\ZZ/2)^2$.

 In this case there are two possibilities for the
component $\mathfrak M_{g;q,[(k'_1, \dots k'_{q-1})]}$, one being the 
divisor of hyperelliptic curves
of genus $g=3$, and the other being the $4$-dimensional locus of 
double covers of elliptic curves.

For $p$ odd, any element of order $2$ in $A' \setminus G$ has exactly 
$6$ fixed points,
and the quotient curve $C_2$ has genus $h_2 = 1 + \frac{p-3}{2}$;
whence we land in the component  $\mathfrak M_{p+1;2,6}$ which has dimension
$ 3 \frac{p-3}{2} + 6 \geq 6 > 3$, and we have a strict inclusion
$$\mathfrak M_{p+1;p,[(0, \dots 0)]} \subset_{ \neq} \mathfrak M_{p+1;2,6}.$$ 

\item
$h_1 = 1, k = 2$ (hence $g=p$) and $G' \cong \ZZ/2$ is generated by a
transformation $ z \mapsto - z + a$;
we have that $[(k_1, \dots k_{p-1})]$ is the class of $k_1 = 1,
k_{p-1} = 1, k_i = 0 \ {\rm for} \ i\neq 1, p-1$;
   $A'$ is a dihedral group $D_p$ and the genus of $C$ equals $p$.
   
   For $p=2$ we have $A' \cong (\ZZ/2)^2$ and  there is only one 
possibility for the
component $\mathfrak M_{g;q,[(k'_1, \dots k'_{q-1})]}$, namely the 
divisor of hyperelliptic curves
of genus $g=3$. 

This  case yields a component of the singular locus 
of $\mathfrak M_3$.

For $p$ odd, any element of order $2$ in $A' \setminus G$ has exactly 
$4$ fixed points,
and the quotient curve $C_2$ has genus $h_2 = 1 + \frac{p-3}{2}$;
whence we land in the component  $\mathfrak M_{p;2,4}$ which has dimension
$  \frac{3}{2} (p-3) + 4 \geq 4 > 3$, and we have again a proper inclusion.

\item
$h_1 = 0, k = 4$,  (hence $g = p-1$)  and $ \ZZ/2 \subset G'$ occurs for arbitrary $p$; if
$\tau \in G'$ is a nontrivial element,
then it is a double transposition of the four branch points. 

If $\tau $ 
permutes $ P_1$ with $P_2$,
and $ P_3$ with $P_4$, then the local monodromies of the four points are
either $m_1= m_2 =1, m_3= m_4 =   p-1$ or $m_1 = 1, m_2 = p-1, m_3 =
n, m_4 = p-n$.

In the first case (where $m_1= m_2 =1, m_3= m_4 =   p-1$), then  $\tau$ lifts and $\tau$ 
and $G$ generate a group
$\cong (\ZZ/2 p)$ for $p$ odd,  in the second case they generate a dihedral group $D_p$.

Thus in the case where the local monodromies are as in the first case (for instance for $p=3$,
since then the two   cases coincide), 
$A' \cong D_p \times( \ZZ/2)$.

In the other case where the local monodromies are (equivalent to) $m_1 = 1, m_2 = p-1, m_3 =
n, m_4 = p-n, \ n \neq 1,  n \neq p-1$ then $A' \cong D_p $.

\item
$h_1 = 0, k = 3$ (hence $ g = \frac{p-1}{2}$ and $ p \geq 5$) and $ \ZZ/2 \subset G'$ occurs for arbitrary $p$, while
$ \ZZ/3 \subset G'$ occurs for $ p \equiv 1 (3)$.

$G' \cong \mathfrak S_3$ cannot occur.

\item
In the case $h_1 = 0, k = 3$, $ \ZZ/2 \cong G'$  we have that $A'$ is cyclic of order $2p$, and contained in the
isotropy subgroup of a unique point of $C$. 

This  is the only case where $A'$ is  contained in the isotropy subgroup of
a  point of $C$.

\end{enumerate}
\end{theo}

\begin{rem}
1) In many cases one can directly show, once one has an inclusion
$\mathfrak M_{g;p,[(k_1, \dots k_{p-1})]} \subset \mathfrak 
M_{g;q,[(k'_1, \dots k'_{q-1})]}$,
  that this inclusion is strict: for instance it suffices that
$ \mathfrak M_{g;q,[(k'_1, \dots k'_{q-1})]}$ does not appear in the 
list given in Theorem
\ref{cornalba}.

2) If the group $A'$ is $ (\ZZ/2 p)$with $p$ odd the component $ 
\mathfrak M_{g;q,[(k'_1, \dots k'_{q-1})]}$
clearly corresponds to the only element of order $2$ in $A'$.
If the group  $A'$ is a dihedral group $D_p$, then we know that all 
elements of order $2$ are conjugated,
hence the corresponding components are of the same topological type.

\end{rem}

\Proof

The quotient curve $C_1$ has a nontrivial group of automorphisms $G'
: = A' /G$, which preserves the set
of $k$ branch points.

Our main observation is however that the curve $C_1$ and the $k$ branch points
are general, therefore
we conclude first of all  that either $h_1 = 2, k = 0$ and $G' \cong
\ZZ/2$ is generated by
the hyperelliptic involution, or $ k \geq 1$.

Since these branch points can be chosen arbitrarily, it follows that
the genus $h_1 \leq 1$. And we can use the result of remark \ref{numbers}.

If $h_1 = 1$, then, since $k \geq 2$, we obtain by the generality assumption on the branch 
points that $k=2$, and
$G'$ has order $2$. In fact, given two general points $x,y $ in a
general elliptic curve, there is
no translation leaving the set $\{x, y\}$ invariant, while there is a unique
transformation $\tau_a : z \mapsto - z + a$ exchanging $x,y$:
the one with $ a = x+y$.

If instead $h_1 = 0$, it must be $k = 4$ (there exists always  permutations
which preserve the cross ratio, given by double transpositions) or $k=3$.

Let's examine now more closely the several possibilities.

{\bf Case 1.}

For $h_1 = 2, k = 0$ we observe that the hyperelliptic involution
acts on the first homology group
as $ - 1$, and $-1$ is an automorphism of $G$, whence the
hyperelliptic involution lifts to an
automorphism of $C$.

 If $ p \neq 2$ we get that $A'$ is a dihedral
group, while for
$p=2$ we have $g=3$ and $ A' \cong (\ZZ/2)^2$. 

In the latter  case we have a
bidouble cover of $\PP^1$
with branch divisors of degrees $0,2,4$, hence $C$ has two more involutions, with respective
quotients of genus $1$, or $0$. In particular $C$ is hyperelliptic 
and the larger component consists either of the
hyperelliptic curves of genus $3$ (this component has dimension $5$), or of the double covers
of elliptic curves (this component has dimension $4$).

In the case where $p$ is an odd prime, all elements in $A' \setminus G$ are conjugate,
and we may take any $\ga \in A' \setminus G$. 

Consider the standard presentation
 $$D_p = \langle x,y ; y^2= x^p = xyxy= 1\rangle, $$
 and observe that all the local monodromies are conjugate of $y$, and denote then 
 by $H$ for each branch point 
 the corresponding subgroup.  So that the fibre over this branch point is in bijection with the set of cosets $ w H$.
 
It is then easy to see that $y$
 has exactly one fixed point lying over each of the $6$ branch points in $\PP^1$:
 since, if $H$  is as above, $ H = z  \{ 1, y\} z^{-1}$,  $$  y w H = w H \Leftrightarrow  w^{-1} y w \in H \Leftrightarrow  w^{-1} y w = z y z^{-1} 
  \Leftrightarrow zw \in  \{ 1, y\}.$$

{\bf Case 2.} 

For $h_1 =1, k =2$, we may assume that the local monodromies of the
two points $x,y$
are equal to $1, p-1$,  since their sum is $0$.
In this case the automorphism $\tau_a$ lifts to the cyclic covering
and we have again a dihedral
group.

For $p=2$ we have again a bidouble cover of $\PP^1$
with branch divisors of degrees $0,2,4$ hence $C$ has two more involutions, with respective
quotients of genus $2$, or $0$. 

If $p$ is odd, the covering $C \ra C' = \PP^1$ is branched in $4$ points, and the normal form of
the monodromies is then $(y,y,y,y,x, x^{-1})$ using the standard presentation of $D_p$
(we do not really need this assertion for the forthcoming argument).
For each involution there is exactly, as before, one fixed point lying above each branch
point of $C_1 \ra  C' = \PP^1$.

{\bf Case 4.}

For $h_1 =0, k =3$, let the local monodromies be $m_1, m_2, m_3 \in
\{1, \dots, p-1\}$. Without loss of
generality we may assume $m_1 = 1$, and recall the obvious inequality    $ 2p-1 \geq 1 +
m_2 + m_3 \equiv 0 (p)$.

Hence $ 1 + m_2 + m_3 = p$ and we may write $m_1 = 1, m_2 = m, m_3 = p -1-m$,
with $ 1 \leq m \leq p-2$.

Let $\tau \in G'$: if $\tau$ has a fixed point, then we assume this
point to be the point $P_1$.
Then $\tau$ transposes $P_2$ and $P_3$, whence it lifts if and only
if $ m = p-1-m$,
i.e., $ m = \frac{p-1}{2}$ (observe that here $p \geq 3$).

If instead $\tau$ cyclically permutes the branch points $P_1 \mapsto
P_2 \mapsto P_3 \mapsto P_1$,
then  $\tau$ lifts if and only if $ m \cdot ( p-1-m) \equiv 1(p)$.

 In the
field $\ZZ / p$ (recall  $p \geq 3$)
this means that $ m^2 + m + 1 \equiv 0$  has a solution, and this occurs if
and only if there exists $m$ which is  is a nontrivial third root of $1$ (since $ (m^2 + m + 1)(m-1) = m^3 -1 $).
 
Thus this holds if and only if   $ p \equiv 1 (3)$.

Observe that the two cases for the existence of such a $\tau$ are mutually exclusive: if the local
monodromies are $ 1,  \frac{p-1}{2},
\frac{p-1}{2}$, then they cannot be all equal, since $p=3$ leads to a
curve $C$ of genus $g=1$.
Therefore it cannot occur that   $G' \cong \mathfrak S_3$.

{\bf Case 3.}

For $h_1 =0, k =4$, let the local monodromies be $m_1, m_2, m_3 , m_4
\in \{1, \dots, p-1\}$.
Without loss of
generality we may assume $m_1 = 1$, and recall that  $ 3p-2 \geq 1 +
m_2 + m_3 + m_4 \equiv 0 (p)$.

Hence $ 1 + m_2 + m_3 + m_4 = p$ or $ 1 + m_2 + m_3 + m_4 = 2 p$ and
we may write $m_1 = 1, m_2 = m,
m_3 =n, m_4 =  p -1-m-n$, or  $m_1 = 1, m_2 = m,
m_3 =n, m_4 =  2p -1-m-n$.

Since the four points have a general cross ratio, $\tau \in G'$ is a
double transposition and
without loss of generality we may assume that $\tau$ permutes $ P_1$
with $P_2$,
and $ P_3$ with $P_4$. 

Whence  $\tau$ lifts if and only if $ m^2 \equiv 1 (p)$,
$ m_4 \equiv m n (p)$. Hence $ m \equiv \pm 1$, and $  1 + m + n + mn
\equiv 0 (p)$,
i.e., $$ ( 1+m) (1 +n) \equiv 0 (p) \Leftrightarrow m \equiv -1 (p)
\ {\rm or } \ n \equiv -1 (p).$$

Hence the solutions are with $ m \equiv -  1$, $n$ arbitrary, $m_4 \equiv - n$,
or $m_1= m_2 =1, m_3= m_4 =   p-1$.

The other assertions follow then easily as before, since if $\tau$ exchanges two points with
the same monodromy, then its lift commutes with $\ga$; else, if it exchanges two
points with opposite monodromies,  then its lift conjugates  $\ga$ with its inverse.

Observe finally that in the second case if it were possible to exchange $P_1$ with $P_3$ it should hold:
$n^2 \equiv 1 (p)$, which is excluded by the assumption $ n \neq 1, n \neq p-1$.

{\bf Assertion 5.}

It is easy to describe explicitly the case where $h_1=0, k=3$ and $\tau$ has order $2$.
Up to an automorphism of the cyclic group, we can assume $m_1 = p-2$,
$m_2=m_3=1$.

In other words, the function field of $C$ is $ \CC (x,z)$, where $z^p = (x^2-1)$.
The involution is $\tau (x) = - x$, which lifts to $C$, and the quotient of
$C$ by $\tau$ is the curve with function field $ \CC (z)$.

Setting $y : = x^2$, $xz : = u$, we have $z^p = (y-1)$,
$ u^{2p} = y^p (y-1)^2$, and $ \CC (x,z) = \CC (y,u)$ expresses
$C$ as a cyclic covering of the $\PP^1$ with function field $ \CC (y)$.

The only fixed point for the cyclic group $G'$ is the one lying over $\infty$,
while the automorphism $\ga$ has three fixed points, lying over
$  x= \infty, x =1, x=-1$. 

$\tau$ exchanges the second and third point of these three,
and leaves fixed other $p$ points, the ones  lying over $x=0$.

Assume now that $A'$ is contained in the isotropy subgroup of a point
$P_0 \in C$, hence in particular $A'$ is cyclic.
From remark \ref{numbers} and the above description of the possible cases, 
we see that we necessarily are in the
case where $h_1=0, k =3$, and $\tau$ leaves one point fixed. Thus we are in case 4,
with $\tau$ of  order $2$.

\qed

   \section{Spaces of cyclic covers of stable curves}

In this section we consider a stable curve $C$ of genus $g \geq 2$
and $\ga \in Aut (C)$ an automorphism,
of order $d$. At a later point we shall make the simplifying
assumption that $d$ is a prime number.

Consider the decomposition of $C$ into irreducible components
$$  C = \cup_{ i \in I}  C_i. $$

We partition the set of irreducible components into three subsets:
$$I_0 : = \{ i \in I | \ga |_{  C_i} = Id_{  C_i} \} ; \  I_1 : = \{
i \in I  \setminus I_0 | \ga ( C_i ) = C_i \} ;  I_2 :
=  I  \setminus \{I_0  \cup I_1\}.$$

Since each irreducible component of the space of such pairs $(C,
\ga)$ will be the union of certain strata, we
try to look immediately for such  strata which are open. For this 
purpose we use
deformation theory.

Recall that the Kuranishi space of such a curve $C$ is smooth and
locally biholomorphic to
$ Ext^1 (\Omega^1_C, \hol_C)$, and the subspace of local deformations
of the pair  $(C, \ga)$ will be
locally biholomorphic to
$ Ext^1 (\Omega^1_C, \hol_C)^G$, where $G$ is the cyclic group
generated by $\ga$.

The local to global spectral sequence yields an exact sequence
$$0 \ra \oplus_i H^1( \Theta_{C_i} ( - \sum _{j \neq i} ( C_i \cap
C_j)) \ra  \Ext^1 (\Omega^1_C, \hol_C)
\ra  (\oplus _{p \in Sing(C)} {\sE xt}^1 (\Omega^1_C,
\hol_C)_p) \ra 0.$$

It is an exact sequence of  $G$-vector spaces, hence the sequence of
$G$-invariants is
also exact.

In particular we have  a surjection:
$$ \Ext^1 (\Omega^1_C, \hol_C)^G \ra  ( \oplus _{p \in Sing(C)} { \sE
xt}^1 (\Omega^1_C,
\hol_C)_p)^G.$$

The first consequence  is that we can smooth all the $G$-fixed nodes
$p$ such that
$$  { \sE xt}^1 (\Omega^1_C, \hol_C)^G_p = { \sE xt}^1 (\Omega^1_C,
\hol_C)_p  \cong \CC.$$

By the well known Cartan's lemma (\cite{cartan}) the action of $G$ 
can be linearized
around $p$,
and, if we set $\zeta := exp ( 2 \pi  \ i / d)$, there are local
holomorphic coordinates
$(x,y)$ such that $C = \{ xy = 0 \}$,  and
either
$$ \ga (x,y) = ( \zeta^{m} x,   \zeta^{n} y)$$
or
$$ \ga (x,y) = ( y ,   \zeta^{2m} x).$$

$ { \sE xt}^1 (\Omega^1_C, \hol_C) \cong \CC $ is identified with the
space of local deformations of the
singularity, i.e., $\CC = \{ t \in \CC \}$ is the parameter space for 
the family of curves
$$\{ (x,y,t) | xy = t \}  .$$

In the first case $G$ acts on the family by
$$ \ga (x,y,t) = ( \zeta^{m} x,   \zeta^{n} y, \zeta^{ m  +n}t).$$

Hence $  { \sE xt}^1 (\Omega^1_C, \hol_C)^G_p = { \sE xt}^1
(\Omega^1_C, \hol_C)_p  \cong \CC$
if and only if $ m + n \equiv 0 (d)$, i.e., exactly when we have a
local family of curves with a $G$-action.

In the second case $G$ acts on the family by
$$ \ga (x,y,t) = ( y   ,  \zeta^{2m} x, \zeta^{2m}t)$$
hence $  { \sE xt}^1 (\Omega^1_C, \hol_C)^G_p = { \sE xt}^1
(\Omega^1_C, \hol_C)_p  \cong \CC$
if and only if $ 2m \equiv 0 (d)$, i.e., exactly when we have a local
family of curves with a $G$-action.

\begin{defin} We shall say that a pair $(C, \ga)$ is  { \bf
simplifiable} if it admits a small deformation to a
pair with a smaller number of nodes, whereas we shall say that $(C,
\ga)$ is { \bf maximal} if it is not
simplifiable.

\end{defin}

\begin{rem} (1) Assume that $i,j \in I_0$ and that $ p \in C_i \cap
C_j $. Then the node $p$ can be
smoothed.

Hence, if $(C, \ga)$ is maximal, $ \forall i,j \in I_0$ we have $
C_i \cap C_j  = \emptyset.$

(2) if $i_0 \in I_0$ and $\exists  p \in C_j \cap C_{i_0}$, then $
\ga (C_j ) = C_j$.

(3) if $p$ is a node such that $\ga (p) = p$, and $ p \in C_i \cap
C_j $, $(i\neq j)$, then either

(3i) $\ga (C_i) = C_i , \  \ga (C_j ) = C_j$ or

(3ii) $\ga (C_i) = C_j$.

In the first subcase, if we set as above $\zeta := exp ( 2 \pi  \ i /
d)$, there are local holomorphic coordinates
$(x,y)$ such that $C = \{ xy = 0 \}$, $C_i = \{ y = 0\}$, $C_j = \{ x
= 0\}$, and
$$ \ga (x,y) = ( \zeta^{m_i} x,   \zeta^{m_j} y).$$

The node is smoothable if and only if $m_i + m_j = d$ (we take $m_i
\in \{ 0,1,\dots , d-1 \}$).

In the second case we have:
$$ \ga (x,y) = ( y,   \zeta^{2m} x),$$ and the node is smoothable if
and only if $ \zeta^{2m}  = 1$.

\end{rem}
\begin{lemma}

If  $d$ is prime, then  each node $p \in C$ not  fixed by $\ga$ can
be smoothed.
Hence, if  $(C, \ga)$ is maximal and $d$ is prime, then every node $p
\in C$ is fixed by $\ga$.

\end{lemma}

\Proof
Since $d$ is prime we see then that the orbit of $p$ has $d$
elements, and there is a bijection
between $G$ and the orbit $ G(p)$.

Look however at the summand of  $ ( \oplus _{p' \in Sing(C)} { \sE
xt}^1 (\Omega^1_C, \hol_C)_{p'})^G$
corresponding to  $ ( \oplus _{ p' \in G (p)} { \sE xt}^1
(\Omega^1_C, \hol_C)_{p'})^G \subset
( \oplus _{ p' \in G (p)} { \sE xt}^1 (\Omega^1_C, \hol_C)_{p'})$.

It is not empty since
   $ \oplus _{ p' \in G (p)} { \sE xt}^1 (\Omega^1_C, \hol_C)_{p'}$
corresponds to the
   representation of $G$ on the orbit $ G(p)$: this yields  a small
deformation smoothing all nodes in the orbit $ G(p)$.

\qed

\begin{prop} Let $\hat{\ga}$ be the permutation of $I$ induced by $\ga$.

If  $(C, \ga)$ is maximal and $d$ is prime, then  $ \hat{\ga }= { \rm
Identity}$
(i.e., $I_2 = \emptyset$).

\end{prop}

\Proof
Observe preliminarly that if  $\hat{\ga}(i) = j \neq i$ and $C_i \cap
C_j \neq  \emptyset$, then $d$
is divisible by $2$. Since then, by the previous lemma, if we take  $ p \in
C_i \cap C_j $,
then $\ga (p) = p$, so that  $\hat{\ga}$ transposes $i$ and $j$.

Since $d$ is prime, $d=2$ and
by the previous remark each node  $ p \in C_i \cap C_j $ is smoothable.

Hence, since we assume $(C, \ga)$ to be maximal,   $C_i \cap C_j  =
\emptyset$.

Let now $ i \in I_2$. There exists, by the connectedness of $C$, a $j 
$ such that $C_i \cap C_j
\neq  \emptyset$.

Let $ p \in C_i \cap C_j$. Since $\ga (p) = p$, $\hat{\ga}$
leaves the set  $\{i,j \} $ invariant.

But $\hat{\ga}(i) = i $ contradicts $ i \in I_2$,while
  $\hat{\ga}(i) = j $ contradicts
  our previous observation. Hence $I_2 = \emptyset$.

\qed

\begin{prop} If $(C, \ga)$ is maximal and $d$ is prime then,  for
each $i \in I_0$, $C_i$ is smooth.

If $ i \in I_1$ and $p $ is a node of $C_i$, then $\ga$ does not
exchange the two branches of $p$.

\end{prop}

\Proof

If $i \in I_0$, then each node $p$ of $C_i$ is smoothable.

If $ i \in I_1$, and $p$ is a node of $C_i$, then we know that $ \ga (p) = p$.

If $\ga$ exchanges the two
branches, then $d = 2$, and as we saw before the node is smoothable.

If $\ga$ does not exchange the two branches, the local action is $
\ga (x,y) = (\zeta^m x , \zeta^n y)$, and
the node is smoothable iff $ m + n = d$.

\qed

We are ready to define the numerical type (it is indeed a combinatorial type, but we call it a numerical type to
keep the analogy with the smooth case) of maximal pairs
$(C,
\ga)$ in the case where $d$ is a prime
number.

\begin{defin} Let $(C, \ga)$ be a maximal pair for $d$ prime.

We attach to $(C, \ga)$ a graph whose vertices correspond to the set
$I$, and whose edges correspond to
the nodes
$p$. Each vertex $i$ has a multilabeling, first of all a labeling by
the genus $g_i$ of $C_i$, and then a
colouring $0, {\rm or }\ 1$, according to
$i \in I_0$, or $ i \in I_1$.

For $i \in I_1$, we associate to $i$ a branching sequence $(k'_1,
\dots, k'_{d-1})$ corresponding to the
fixed points of $\ga | _{C_i}$ which are not nodes.

For each edge $p$ connecting $i$ and $j$, $ i \in I_1$, $ i \neq j$,
we give  labels $m (p,i)  \in \{ 1,\dots ,
d-1 \}$,
$m (p,j)  \in \{ 0,1,\dots , d-1 \}$ according to the local action at
the fixed point $p \in C_i$, respectively
$p \in C_j$.

If instead $i=j$, we look at the action on the two branches and
obtain an unordered pair  $n_1(p,i), n_2(p,i)$.

\end{defin}

In order to understand the notion of admissibility which will be given next,
observe that for each $j \in I_1$ we shall consider a curve $C''_j$
which is the
normalization of $C_j$, hence the genus $C''_j$ equals $g_j$ plus the
number of nodes of $C_j$.
Then $G$ acts on $C''_j$ and we denote by
$C'_j : = C''_j / G$ the quotient curve,  by
$g'_j$ the genus of $C'_j$, and by
   $r_j$ the number of branch points of $C''_j \ra C'_j$.

\begin{rem}
$C_j$ and $C_j / G$ can be easily reconstructed by the marking of
certain pairs of branch points
on $C'_j$.
\end{rem}

\begin{defin} Let  $d$ be a prime number.

Then an admissible automorphism graph is  a connected graph with the
following properties.

It has set of vertices $ I = I_0 \cup I_1$, and  each vertex $ i \in
I_0$ is labelled by an integer $g_i$, while
each vertex $ i \in I_1$ is labelled by an integer $g_i$, and by a
branching sequence $(k'_1(i), \dots,
k'_{d-1}(i))$.

No edge can connect two vertices in $I_0$, but an edge can connect a
vertex $i \in I_1$ with itself (i.e., the
graph has loops).

Each edge $p$ connecting $i$ and $j$, $ i \neq j$, is labelled by $m
(p,i)  \in \{ 0,1,\dots , d-1 \}$, and $m
(p,j)  \in \{ 0,1,\dots , d-1 \}$ in such a way that  $m (p,i) = 0$
if and only if $i \in I_0$, and moreover $
m(p,i) + m (p,j) \neq d$.

If instead $i=j$, we label the loop by  an unordered pair  of 
integers in  $ \{ 1,\dots , d-1 \}$
$$n_1(p,i),
n_2(p,i) \neq 0$$ such that
   $n_1(p,i) + n_2(p,i) \neq d$.

   Define, for each $i \in I_1$,  the  branching integer $k_m(i)$ as the
sum of $k'_m(i)$
   with the number $k''_m(i)$ of occurrences of the integer $m$ among the
    integers $m(p,i)$, for $p$ an edge
connecting $i$ with $j \neq i$,
   or among the pairs of integers $n_1(p,i), n_2(p,i) $, for $p$ a 
loop based at $i$.

   Then the sequence $(k_1(i), \dots, k_{d-1}(i))$ must be {\bf admissible}
in the sense that there
   exists a non negative integer $g'_i$ such that, setting
$g''_i := g_i + \nu_i$. $\nu_i$ being the number of loops
based at the vertex $i$, and setting $ k(i) : = \sum_m k_m(i)$, we have

   $$ 2(g''_i-1) = d  \{  2 (g'_i-1) +  k (i) (1 - \frac{ 1}{d}) \}.$$

   The genus $g$ of the graph is as usual defined as
   $$ g : = \sum_{i \in I} g_i + b^1,  $$
   where $b^1$ is the first Betti number of the connected graph.

\end{defin}

\begin{rem} There is an obvious action of $( \ZZ/d)^*$ on the
branching sequences and on the edge
labelings, and the equivalence classes of the admissible graphs for
this action are denoted the { \bf
numerical types} of the maximal pairs
$(C,G)$ where $G$ is a cyclic group of automorphisms of prime order $d$.

\end{rem}

\begin{theo}\label{stableirred} The pairs $(C,G)$ where $C$ is a
stable projective curve of genus $g \geq
2$, and $G$ is a finite cyclic group of prime order
$d$ acting faithfully on $C$ with a given numerical type associated
to an admissible automorphism graph
$\sG$
   are  parametrized by a  non empty connected complex manifold
$\sT_{g; d, [\sG] }$.

The image  $\overline{\mathfrak M}_{g; d, [\sG]  }$ of $\sT_{g; d,
[\sG] }$ inside the compactified moduli
space
   $\overline{\mathfrak M_g}$ is a locally closed subset of the same
dimension whose closure consists of the
coverings whose numerical type can be simplified to the numerical
type of $\sT_{g; d, [\sG] }$.

If $\sT_{g; d, [\sG] }$ contains only stable singular curves, then
$\overline{\mathfrak M_{g; d, [\sG]  }}$ is
not a divisor in the moduli space  $\overline{\mathfrak M_g}$, unless
we are in the following two cases:
\begin{enumerate}
\item
$d=2$, $ C = C_1 \cup C_2$, where $1 \in I_0$, $2 \in I_1$, and
$g_2=1$ (elliptic tail).
\item
$d=2$, $ C = C_1 \cup C_2$, where $1,2 \in I_1$, and $g_i=1$ ($g=2$ case).

\end{enumerate}
\end{theo}

\Proof

We consider $\sT_{g; d, [\sG] }$ as a product of  two products of
Teichm\"uller spaces: firstly
$$\Pi_{i \in I_0}  \sT_{g_i, r_i} ,$$ where $r_i$ is the number of
edges which touch
the   vertex $i \in I_0$, and
secondly
$$\Pi_{j \in I_1}  \sT_{g'_j, r_j} .$$

In the second case, for each $j \in I_1$, and for each point in $ 
\sT_{g'_j, r_j}$ we
construct a curve $C''_j$
which has an automorphism $\ga$ of order $d$ with quotient a curve $C'_j$
of genus $g'_j$, and such that $C''_j  \ra C'_j $ is branched on
$r_j$ points, and with the
branching indices determined by $\sG$.

Then we construct the family of curves $C_j$ from the family of curves  $C''_j$
glueing  certain pairs of ramification points according to the pattern
determined by $\sG$.

Observe that the family of such coverings $C''_j \ra C'_j$ is an
irreducible family in view of Theorem
\ref{irreducible}, since the numerical type determines the
branching datum, and the local  monodromy at the
points of $C'_j$ corresponding to the nodes of $C_j$ is also
determined by the admissible automorphism
graph.

Hence the same holds for the family of such coverings $C_j$ (for $j
\in I_1$), and we finally
obtain in this way a family with connected base of curves $C$ with an
automorphism $\ga$ of topological
type determined by the graph $\sG$.

Since by our assumption for each curve $C$ in the family we have
$$( \oplus _{p \in Sing(C)} { \sE xt}^1 (\Omega^1_C, \hol_C)_p)^G = 0,$$
   it is easy to see that our family
equals the Kuranishi family $ \Ext^1 (\Omega^1_C, \hol_C)^G $ of
pairs $(C, \ga)$ at each point.

Hence the image  $\overline{\mathfrak M_{g; d, [\sG]  }}$ of $\sT_{g;
d, [\sG] }$ is a locally closed subset of
the same dimension as $\sT_{g; d, [\sG] }$. Its closure consists of
pairs $(C,G)$ corresponding to pairs
$(C, \ga)$ which admit a $G$-invariant local deformation containing a
maximal pair $(C, \ga)$.

If $\sT_{g; d, [\sG] }$ contains only stable singular curves, then
the locus $\overline{\mathfrak M_{g; d,
[\sG]  }}$ is not a divisor  if the general curve has at least two nodes.

If  $C$ is stable and has only one node, then if $C$ is irreducible
its normalization has no automorphisms
provided $C$ is general. If instead $ C = C_1 \cup C_2$, we should
have that all smooth curves $C_1$ of
genus $g_1$ and  all smooth curves $C_2$ of genus $g_2$ occur.

Assume that $1 \in I_0$ and $2 \in I_1$: then it must be $g_2 = 2$ or
$g_2 = 1$. Since however the node
must be, for $g_2=2$, a fixed point for the hyperelliptic involution,
we do not have a divisor for $g_2 =2$
and this case must be excluded.

If instead $ 1,2 \in I_1$, both curves would have genus equal to $2$
or $1$. If however $g_1 = 2$, then the
node must be a fixed point for the hyperelliptic involution, and we
no longer have a divisor.

\qed

We need the following result for the forthcoming theorem.

\begin{lemma}\label{tails} Let $C$ be a stable curve with elliptic tails $E_1,
\dots E_r$.  Then the automorphism group
of $C$ is a direct sum $ (\oplus _{i=1, \dots r} \ZZ/m_i)^r \oplus
\Ga$ where the first addendum is
generated by the multiplication by an $m_i$-th root of unity around
the node $p_i$ of each elliptic tail
$E_i$ (hence $m_i \in \{2,4,6\}$). The quotient of the Kuranishi
family of $C$ by $Aut(C)$  is then
singular unless $\Ga$ is trivial and $m_i =2 \ \forall i$.

\end{lemma}

\Proof It suffices to define $\Ga$ as the subgroup which acts as the
identity on each elliptic tail $E_i$.

Recall that, by the cited theorem of Chevalley (\cite{chevalley}), 
the quotient of a smooth
manifold of dimension $N$ by a finite group
$\Ga$ acting with a fixed point $P$ is smooth at the image point of $ 
P$  if and only if the group is generated by
pseudoreflections (these are transformations which
are biholomorphic to linear maps with exactly $N-1$ eigenvalues equal to $1$).
In particular, if the quotient is smooth the
locus of fixed points is a union of divisors.

The only pseudoreflections correspond to elliptic tails and to the
case where the automorphism of $E_i$ is of
order 2. The assertion follows now immediately.

\qed

\begin{defin}\label{enlargement} {\bf (Enlargement)}.

Let $\ga$ be  an automorphism of prime order $d > 1$, 
with  $ (C, \ga)$  maximal (hence $ I = I_0 \cup I_1$).

{\bf Enlargement of type  1).}

Assume there is a
component $C_j$, where $j \in I_1$, such that $C_j$ does not intersect
any component $C_i$ with $i \in
I_0$. Then we can consider an automorphism
$\ga''$ such that $\ga | _{C_j} = {\rm Identity}$, and $\ga'' = \ga$
on the other components.

If $C_j$ has nodes, we obtain a non maximal pair $(C, \ga'')$, where
   $\ga''$ has the same order of $\ga$, but if we smooth the nodes
of $C_j$, we obtain a maximal pair $(C'', \ga'')$.
Denote by $\sG''$ the associated graph.

Then  the closed   subvariety $\overline{\mathfrak M_{g; d, [\sG]
}}$ is  contained in  the closed   subvariety $\overline{\mathfrak
M_{g; d, [\sG'']  }}$,
and properly contained unless
$C_j$ is an elliptic tail and $d=2$ (if $C_j$  is smooth of genus $0$,
observe that it intersects the other components in at least 3 points, which are
fixed by $\ga$, contradicting $j \in I_1$).
\smallskip

{\bf Enlargement of type  2).}

Assume that  there is a
component $C_j$ with $j \in I_1$, such that $C_j$  intersects some
components $C_i$ with $i \in
I_0$. Assume further that $I_1 \neq \{j\}$.

Then we can consider an automorphism
$\ga''$ such that $\ga | _{C_j} = {\rm Identity}$, and $\ga'' = \ga$
on the other components.
We obtain a non maximal pair $(C, \ga'')$, where
   $\ga''$ has the same order of $\ga$.  If we first smooth the nodes of $C_j$,
we obtain a non maximal pair $(C'', \ga'')$; however, if we smooth the
union of $C''_j$
with the  components $C_i$ with $i \in
I_0$ and  which intersect $C''_j$, we obtain a maximal pair $(C''', \ga''')$.
Denote by $\sG'''$ the associated graph.

Then  the closed   subvariety $\overline{\mathfrak M_{g; d, [\sG]
}}$ is  contained in  the closed   subvariety $\overline{\mathfrak
M_{g; d, [\sG''']  }}$,
and properly contained unless
$C_j$ is an elliptic tail and $d=2$.

{\bf Maximal enlargement.}

Take a
component $C_j$ with $j \in I_1$, assume further that $I_1 \neq \{j\}$.

Then we can consider an automorphism
$\ga''$ such that $\ga | _{C_i} = {\rm Identity}$ for $i \neq j$, and $\ga'' = \ga$
on $C_j$.

We obtain a non maximal pair $(C, \ga'')$, where
   $\ga''$ has the same order of $\ga$.  If we  smooth the nodes of $C \setminus C_j$,
we obtain a maximal pair $(C'', \ga'')$.

Denote by $\sG'''$ the associated graph.

Then  the closed   subvariety $\overline{\mathfrak M_{g; d, [\sG]
}}$ is  contained in  the closed   subvariety $\overline{\mathfrak
M_{g; d, [\sG''']  }}$,
and properly contained unless $d=2$ and, $\forall h \in I_1 \setminus \{j\}$,
$C_h$ is an elliptic tail.

\end{defin}

Clearly a maximal enlargement can be obtained as a sequence of enlargements of
type 1 and 2.

Before stating the main theorem of this section, let's discuss two cases
to which we shall refer  as to ` the exceptional cases'.

\begin{defin}
Consider a curve $C$ as in 5 of theorem \ref{cornalba}, i.e., with a cyclic group $A'$
of automorphisms of order $2p$, where $p \geq 3$ is a prime number.

Denote by $\ga'$ an automorphism in $A'$ of order $p$, and by $\ga$ the automorphism of order
$2$.

Denote by $P_0$ the only fixed point of $A'$ on $C$, and let $P_1, P_2$ be the other
two fixed points of $\ga'$, which are exchanged by $\ga$.

{\bf II-a)}. 

Define $C'_0$ to be the nodal curve obtained by $C$ identifying $P_1, P_2$,
and let $C'_3$ be another smooth curve of genus at least $1$, and $P_3 \in C'_3$ an arbitrary point.

Let $C'$ be the stable curve obtained as $ C' = C'_0 \cup C'_3$, identifying  $ P \in  C'_0 $
with $ P_3 \in  C'_3$, and let $\ga'$ be the automorphism induced by $\ga'$ 
on $C'_0$, extended as the identity on $C'_3$.

Let $\ga $ be the automorphism induced by  $\ga $ on $C'_0$, extended as the identity on $C'_3$.

Then we set $(C'', \ga)$ to be a smoothing of $C'$ at the node of $C'_0$. 

Here, $C''_0$ is smooth 
hyperelliptic of genus $\frac{p+1}{2}$.

{\bf II-b)}. 

Define $C'_0$ to be the curve  $C$, let $C'_3$ be another smooth curve of genus at least $1$, and $P_3 \in C'_3$ an
arbitrary point, and  consider moreover two isomorphic 1-pointed smooth curves of genus at least 1
  $(C'_1, P'_1) \cong (C'_2, P'_2) $.

Let $C'$ be the stable curve obtained as $ C' = C'_0 \cup C'_1 \cup C'_2 \cup C'_3$, obtained identifying  $ P \in 
C'_0 $ with $ P_3 \in  C'_3$, and $P_h \in C'_0 $ with $ P'_h \in  C'_h$, for $ h= 1,2$.

Let $\ga'$ be the automorphism induced by $\ga'$ 
on $C'_0$, extended as the identity on the other components $C'_i$.

Let $\ga $ be the automorphism induced by  $\ga $ on $C'_0$, extended as the identity on $C'_3$,
and exchanging  $C'_1$ with $C'_2$ according to the given isomorphism and its inverse.

Then we set $(C'', \ga)$ to be a smoothing of $C'$ at the nodes corresponding to $P_2, P_3$.

\end{defin}

\begin{theo}\label{stablesing} Assume that $g \geq 2$, and consider the
   closed   subvarieties $\overline{\mathfrak M_{g; d, [\sG]  }}$
inside the compactified moduli space
   $\overline{\mathfrak M_g}$, such that

   \begin{enumerate}
   \item
   $d$ is a prime number
    \item
    the cyclic group $G$ either has order $ d \neq 2$ or it acts
   trivially on the elliptic tails.
\item the subset $I_1$ contains exactly one element
\item
$I_0$ is not empty
(hence   $\overline{\mathfrak M_{g; d, [\sG]  }}$ contains
only singular stable curves).
\item
$ \mathfrak M_{g; d, [\sG]  } $ is not one of the two exceptional cases II-a, II-b 
for $ (C', \ga')$.

   \end{enumerate}

The above   components $\overline{\mathfrak M_{g; d, [\sG]  }}$   are
then all distinct, for different $d$
and different topological types, and provide the irreducible
components of $ Sing (\overline{\mathfrak
M_g)}$ which do not intersect $\mathfrak M_g$.

\end{theo}

\Proof
$\overline{\mathfrak M_g}$ is locally the quotient of the Kuranishi
family of a stable curve $C$ by the
group $Aut (C)$.

Hence $\overline{\mathfrak M_g}$ is smooth unless $C$ has an
automorphism of prime order.

Moreover, as we already recalled, by Chevalley's theorem, the
quotient  is smooth if and only if the action of the group $Aut(C)$ 
is generated by pseudoreflections.

By our previous remark the only
pseudoreflections correspond to reflections on an elliptic tail. Hence, by
  lemma \ref{tails}, it suffices to take care of the cases where $d 
\neq 2$, or, if $d=2$, we can assume that
$\ga$ acts trivially on the elliptic tails.

By the definition of maximal enlargement given in  \ref{enlargement} we may restrict ourselves to consider
only irreducible components
satisfying the assumption that the subset $I_1$ contains exactly one element.

We want first to see when two  irreducible components  $\overline{\mathfrak
M_{g; d, [\sG]  }}$ and
$\overline{\mathfrak M_{g; d', [\sG']  }}$,
satisfying our assumptions,  can be contained into each other.

We assume   $\overline{\mathfrak M_{g; d',
[\sG']  }} \subset \overline{\mathfrak M_{g; d, [\sG]  }}$ .

Let $C'$ be a general curve in $\mathfrak M_{g; d', [\sG']  }$, and
$C$ a general curve in $\mathfrak M_{g; d,
[\sG]  }$. Since  $C$ is a smoothing of $C'$, we have an automorphism
$\ga$ of $C'$ such that
the pair $(C', \ga)$ deforms to the pair $(C, \ga)$.

Observe that there is a unique component $C'_j$ with $ j \in I_1'$.
Since $C'$ is general,
and $(C', \ga')$ is maximal, $C'_j$ intersects all other components
$C'_i$, while each of these other components $C'_i$ is smooth,
and intersects only the component
$C'_j$.

We infer easily from the above observation that there are only  two 
possible cases.

Case a) : $\ga (C'_j) = C'_j$

Case b) : $C'$ has  exactly two components, $C'_j$ and $C'_i$, and
$\ga$ exchanges them.

Case b) leads to a contradiction, since then $d=2$ and $C$ would be smooth
(while $I_0 \neq \emptyset$).

In case a), if $\ga$ were the identity on $C'_j$, we would derive
a contradiction.
In fact, for  each component $C'_i$ with $i \in I_0$ such that $\ga$ is not the identity
on $C'_i$, we have  that $C'_i$ is smooth, and all the points of intersection with
$C'_j$ can be chosen freely. This implies that either $C'_i$ is an elliptic tail, and $\ga$
the elliptic involution, or  $C'_i \cong \PP^1$, and $\ga$ has at least three fixed points
on  $C'_i$ . The second alternative  implies that $\ga$ is the identity of  $C'_i$,
a contradiction. Hence, for all such $ i \in I_0$, the first alternative holds,
contradicting our assumption 2.

    We conclude that both
automorphisms $\ga$,
$\ga'$ leave then  $C'_j$ invariant and are different from the identity.

We have that $C$ contains an irreducible component $C_i$ on which
$\ga$ is the identity.
It follows that $C_i$ comes from deforming a component $C'_i$ with $
i \in I'_0$.

Then both  $\ga$,
$\ga'$ leave $C'_i$ pointwise fixed, in particular a node $P \in C'_i
\cap C'_j$.

Since  $\ga$,
$\ga'$ belong to the isotropy group of $P$, a cyclic subgroup,
 either 

I) $\ga = \ga'$ on $C'_j$, so that 
in particular $ d = d'$, or

II) we are  in the exceptional situation $h_1=0, k=3$ of Theorem  \ref{cornalba}, and 
 $\ga$ has order $d=2$.

If I) holds, then it follows  that, for each component $C'_i$ with $ i \in I'_0$,
$C'_i$ is  left invariant by $\ga$, and there are two possibilities:

A)  $\ga$ is  the identity on $C'_i$ and no nodes $P \in C'_i \cap
C'_j$ are  smoothed;

B) $\ga$ is not the identity on $C'_i$ and some nodes $P \in C'_i \cap
C'_j$ are smoothed;

Case B) leads however to a contradiction, by the generality of $C'$, because
$C'_i$ admits no nontrivial automorphism unless it is an elliptic
tail and $d'=2$.
But this possibility is excluded by our assumptions.

The conclusion is that $C= C'$ and $\ga = \ga'$, as we wanted to show.

Let's consider now the second possibility II).

Here, as we saw in the proof of theorem \ref{cornalba},
$\ga'$ has  exactly 3 fixed points on the normalization of $C'_j$,
one being our $P$, fixed by $\ga$ also, and the other two, $P_1, P_2$ being exchanged by $\ga$.

The first conclusion is that there is only one component $C'_i$ sent to itself by $\ga$.

Since  moreover, for each node of $C'_j$, $\ga'$ does not exchange the two branches,
the only possibilities are :

II-a) $C'_j$ has a node, and  there is only another component $C'_i$, meeting
$C'_j$ precisely in one point $P$. 

II-b) $C'_j$ is smooth, and  there are just three other components: $C'_i$ meeting
$C'_j$ precisely in one point $P$, and then 
$C'_{i_1}, C'_{i_2}$, exchanged by $\ga$ and such that  $C'_j \cap C'_{i_h}$ consists precisely of
 the point $P_h$. 

The above possibilities correspond exactly to the two exceptional cases for $(C', \ga')$, 
hence our proof is finished.

\qed

\section{ Some open questions.}

The first natural question that would be worth to investigate is: 
what is the description
of the numerical type of an automorphism of a stable curve $C$  of 
non prime order $d$?

We have seen in section 4 that the assumption that $d$ be a prime 
number is repeatedly
used, so the combinatorial description is likely to be rather more complicated.

Morally, however, one should expect that again  a similar  result to
theorem \ref{stableirred} holds true in the non prime case.

More generally, an interesting problem is the  investigation of the group of automorphisms
of a stable singular curve.

\begin{rem}
Consider a stable curve consisting of a rational smooth component $C_0$,
intersecting $g$ elliptic tails $C_1, \dots C_g$ in nodes $p_1,  \dots p_g $.

Then, if $C_1, \dots C_g$ and  $p_1,  \dots p_g \in C_0$ are general,
$ Aut (C)$ has cardinality at least $ 2^g$. If instead the elliptic curves
are equianharmonic, and the points $p_1,  \dots p_g $. are roots of 
unity in the complex line
$\CC$, then $ Aut (C)$ has cardinality  $ (2g) \cdot 6^g$.

This number is by far larger than the Hurwitz bound $84 (g-1)$ for the
cardinality  of  $ Aut (C)$  for a smooth
curve of genus $g$.

Concerning the questions about determining  the Hurwitz bound for stable curves,
and  the geometrical
description of the stable curves $C$ of genus $g$ such that $ Aut (C)$
attains the maximal allowed cardinality, which we had posed in the first version of this paper, 
we have been informed by Gavril Farkas that these issues have  been
thoroughly investigated and completely solved by van Opstall and Veliche (see \cite{VO1}, \cite{VO2},
\cite{VO3}).

\end{rem}

\bigskip

\noindent
{\bf Acknowledgements.}
I  would like to  thank Barbara Fantechi, who posed the problem of proving
Theorem \ref{irreducible} during her  talk at MSRI in February 2009.
I regret that, after our initial discussion at MSRI, she 
could not find the time
to collaborate together on the ensuing research project.

I would also like to thank Fabio Perroni for some useful remark.


\end{document}